\documentclass[11pt]{article}
\usepackage{amsmath}
\usepackage{mathrsfs}
\usepackage{amsfonts}

\usepackage{amsfonts, amsmath, amssymb}
\usepackage{amssymb,amsfonts,amsmath,color,
latexsym, epsfig,cite, psfrag,eepic,colordvi}
\usepackage{amscd,graphics}
\usepackage{comment}
\newcommand{\qed}{\hfill \rule{4pt}{7pt}}
\newcommand{\pf}{\noindent {\bf Proof.} }

\newtheorem{theorem}{Theorem}[section]
\newtheorem{lemma}[theorem]{Lemma}
\newtheorem{coro}[theorem]{Corollary}

\newtheorem{conjecture}[theorem]{Conjecture}

\begin{document}

\title{A note on exact  minimum degree threshold for fractional perfect matchings}

\author{Hongliang Lu \footnote{luhongliang@mail.xjtu.edu.cn; partially supported by the National Natural
Science Foundation of China under grant No.11871391 and
Fundamental Research Funds for the Central Universities}
\\School of Mathematics and Statistics\\
Xi'an Jiaotong University\\
Xi'an, Shaanxi 710049, China\\
\smallskip\\
Xingxing Yu \footnote{yu@math.gatech.edu; partially supported by NSF grant DMS-1954134}\\
School of Mathematics\\
Georgia Institute of Technology\\
Atlanta, GA 30332, USA\\}

\date{}

\maketitle

\begin{abstract}
R\"odl, Ruci\'nski, and Szemer\'edi determined the minimum
$(k-1)$-degree threshold for the existence of fractional perfect matchings in
$k$-uniform hypergrahs, and
 K\"uhn, Osthus, and Townsend extended this result by asymptotically
 determining the $d$-degree threshold for the range $k-1>d\ge k/2$.
In this note, we prove the following exact degree threshold:
Let $k,d$ be positive integers with $k\ge 4$ and $k-1>d\geq k/2$, and
let $n$ be any integer with $n\ge
k^2$. Then any  $n$-vertex $k$-uniform hypergraph with minimum $d$-degree $\delta_d(H)>{n-d\choose k-d} -{n-d-(\lceil n/k\rceil-1)\choose
  k-d}$ contains a fractional perfect  matching.  This lower
bound on the minimum $d$-degree is best possible. We also determine optimal minimum $d$-degree conditions which guarantees the existence of fractional matchings of size $s$, where $0<s\le n/k$ (when $k/2\le d\le k-1$), or
with $s$ large enough and $s\le
n/k$ (when $2k/5<d<k/2$).

\end{abstract}

\section{Introduction}
Let $k$ be a positive integer and let  $[k]:=\{1, \ldots, k\}$.  For a set
$S$, let ${S\choose k}:=\{T\subseteq S: |T|=k\}$. A {\it hypergraph} $H$
consists of a vertex set $V(H)$ and an edge set $E(H)$ whose members
are subsets of $V(H)$, and $H$ is said to be {\it $k$-uniform} if $E(H)\subseteq {V(H)\choose k}$. A $k$-uniform hypergraph is also
called a {\it $k$-graph}.
A {\it matching} in a hypergraph $H$ is a set of pairwise disjoint
edges of $H$, and a matching in $H$ is \emph{perfect} if the union of all edges in
the matching is $V(H)$. We use $\nu(H)$ to
denote the largest size of a matching in $H$. A {\it maximum
matching} in  $H$ is a matching in $H$ of size $\nu(H)$. %For positive integer $n$, let $[n]=\{1,\ldots,n\}$

There has been  much activity on degree thresholds for matchings of
certain size in uniform hypergraphs. Let $H$ be a hypergraph.
For  $S\subseteq V(H)$,  let $N_{H}(S)=\{T\subseteq V(H)\setminus S\ :\ T\cup S\in
E(H)\}$ and let $d_H(S):=|N_H(S)|$.
For integer $d\ge 0$, let
$\delta_d(H) =\min \left\{d_H(S):
S\in {V(H)\choose d}\right\}$, which is  the {\it minimum $d$-degree}  of
$H$. Then $\delta_0(H)=e(H)$, the number of edges in $H$.
%$\delta_1(H)$ is the minimum {\it vertex degree} of $H$, and
%$\delta_{k-1}(H)$ is the minimum {\it codegree} of $H$.
 For integers $n, k, d,s$ satisfying $0\leq d\leq k-1$ and $0<s\leq n/k$,
let $m_d^s(k, n)$ denote the minimum integer $m$ such that every $k$-graph $H$ on $n$ vertices with $\delta_d(H)\geq  m$
has a matching of size $s$.

R\"odl, Rucin\'ski, and Szemer\'edi \cite{RRS09}
determined  $m_{k-1}^{n/k}(k,n)$ for all integers
$k\ge 3$ and $n\in k\mathbb{Z}$ sufficiently large.
 Given positive integers $k,d$ with  $k\geq 4$ and $k-2\geq d\geq k/2$,
Treglown and Zhao \cite{TZ12,TZ13} showed that
$m_d^{n/k}(k,n)\sim\frac{1}{2}{n-d\choose k-d}$.

One approach to finding a large matching in a $k$-graph is to first find a
large fractional matching in the $k$-graph, and then convert that
fractional matching to a matching. This approach has been used quite
often, for example, in \cite{AFHRS12,Han16,KOT14}.
A \emph{fractional matching} in a $k$-graph $H$ is a function $f : E(H) \rightarrow [0, 1]$ such that, for each $v\in V(H)$,
$\sum_{\{e\in E(H)  : v\in e\}} f(e) \leq 1$. The {\it size} of $f$ is
$\sum_{e\in E(H)} f(e)$, and
$f$ is a \emph{fractional perfect matching} if it has size
$|V(H)|/k$. We use $\nu'(H)$ to denote the maximum size of a
fractional matching in $H$.
For integers $n, k, d$ and positive rational number $s$ satisfying $0\leq d\leq k-1$ and $ s\leq n/k$,
let $f_d^s(k, n)$ denote
the minimum integer $m$ such that every $k$-graph $H$ on $n$ vertices with $\delta_d(H)\geq  m$
has a fractional  matching of size $s$.

 Alon  et al. \cite{AFHRS12} provided a connection between the
 parameters $m_d^{s}(k, n)$ and $f_d^s(k, n)$.
Let $k,d$ be integers such that $1\leq d\leq k -1$ and let $n$ be a
sufficiently large integer. If there exists $c^*>0$ such that
$f_d^{n/k}(k, n)\sim c^*{n-d\choose k-d}$, then $m_d^{n/k}(k, n)\sim
\max\{c^*, 1/2\}{n-d\choose k-d}$. (For
 integer-valued functions $h_1(n),h_2(n)$,  we write $h_1(n)\sim
 h_2(n)$ if  $\lim_{n\rightarrow \infty} h_1(n)/h_2(n)=1$.)
In the same paper, they show a way to convert a large fractional
matching to a matching using an absorbing technique and a two-round
randomization technique; while K\"{u}hn, Osthus, and
Townsend \cite{KOT14}  used the weak regularity lemma for hypergraphs to
show  $m_d^{an}\sim (1-(1-a)^{k-d}){n-d\choose k-d}$,  where $0\leq a<\min\{(k-d)/2,(1-\varepsilon)n/k\}$ and $\varepsilon>0$ is a constant.

R\"odl, Ruci\'nski, and Szemer\'edi \cite{RRE06} proved   that
$f_{k-1}^{n/k}(k,n)=\lceil n/k\rceil$, which is much smaller than
$m_{k-1}^{n/k}(k,n)$ when $n\in k\mathbb{Z}$ (which is approximately $n/2$).
K\"uhn, Osthus, and Townsend \cite{KOT14} determined $f_d^s(k,n)$
asymptotically when $s\le n/(2(k-2))$ or $d\ge k/2$.

Alon  et al. \cite{AFHRS12} conjectured that for all $1\le d\le
k-1$, $f_d^{n/k}(k,n)\sim (1-(1-1/k)^{k-d}){n-d\choose k-d}$,  and
proved it for $k\ge 3$ and $k-4\le d\le k-1$.
In this note, we determine the exact value of $f_d^{n/k}(k,n)$ for certain ranges of $d$,
using a result of Frankl \cite{Fr13}  and a result of  Frankl and
Kupavskii \cite{FK18}. This is a special case of the following
result.

\begin{theorem}\label{main-frac}
Let $k, d$ be integers such that  $k\geq 4$ and  $2k/5<d\le k-1$. There exists  $s_0=s_0(k,d)$ such that, for  any integer $n$ with $n\ge 2k^2$
and every rational $s$ with  $s_0<s \le n/k$,  $f_d^{s}(k,n)=
{n-d\choose k-d} -{n-d-(\lceil s\rceil-1)\choose k-d}+1$. Actually, we can take $s_0(k,d)=0$  when $d\ge k/2$ and $s_0(k,d)=1$ when $2k/5<d<k/2$.
\end{theorem}

In Section 2, we prove a technical result, Lemma~\ref{Erdos-frac},  about fractional matchings.  In Section 3, we give a short proof of
Theorem~\ref{main-frac} by applying Lemma~\ref{Erdos-frac},  a result of Frankl
(Lemma~\ref{Fran13}), and a result of
Frankl and Kupavskii (Lemma~\ref{FrKu18}).
We will also discuss other related work on asymptotic and
exact bounds for $f_d^{s}(k,n)$ in Section 4.

\section{Fractional matchings}

One of the ideas in our proof is to use the strong duality between the size of a largest
fractional matching in a hypergraph and the size of a smallest fractional vertex cover
of that hypergraph. This idea has been already explored before,  e.g., see
\cite{AFHRS12,KOT14}.
Let $H$ be a hypergraph. A \emph{fractional vertex cover} of $H$ is a function $\omega : V(H) \rightarrow [0, 1]$, such that for
each $e\in E(H)$ we have $\sum_
{\{v\ :\  v\in e\}} \omega(v) \geq 1$. The {\it size} of $\omega$ is
$\sum_{v\in V(H)} \omega(v)$.  We use $\mu(H)$ to denote the minimum size of a
fractional vertex cover in $H$.
Note that $\nu'(H)=\mu(H)$ for any hypergraph $H$, as they are optimal solutions of
two dual linear programs. In our proof of Theorem~\ref{main-frac}, we
will use this fact to  transform the
fractional matching problem on $H$ to one on another hypergraph $H'$.

First, observe that ${n-d\choose k-d}-{n-d-(\lceil s\rceil-1)\choose
  k-d}+1$ is a  lower bound for $f_d^s(k,n)$. For convenience, we
state it as lemma below. The construction involved in the proof is
standard, e.g., see equations (3) and (4) in \cite{AFHRS12}.

\begin{lemma}\label{low-bound}
Let $k,d$ be  integers such that $k\geq 2$ and $0\leq d\leq k-1$. Then,
for any integer $n$ with $n\ge k$ and any rational $s$ with  $0< s \le n/k$, $f_d^s(k,n)\geq {n-d\choose k-d}-{n-d-(\lceil s\rceil-1)\choose k-d}+1$.
\end{lemma}

\pf Let $H_k(n,s)$ be the $k$-uniform hypergraph  with vertex set
$[n]$ and edge set consisting of all $k$-element subsets of $[n]$
which have non-empty intersection with the  subset $[\lceil
s\rceil-1]$.

First, suppose $0<s\le 1$. Then, by definition, $H_k(n,s)$ has no
edge and, thus, has no
fractional matching of any positive size.  Therefore, in this case,
$f_d^s(k,n)\geq 1={n-d\choose k-d}-{n-d-(\lceil s\rceil-1)\choose k-d}+1$.

Hence, we may assume $s>1$. Then
\begin{align*}
\delta_d(H_k(n,s))={n-d\choose k-d}-{n-d-(\lceil s \rceil-1)\choose k-d}.
\end{align*}
 Let $\omega:[n]\rightarrow [0,1]$ such that $\omega(x)=1$ for all $x\in [\lceil s\rceil-1]$ and $\omega(x)=0$ for all $x\in [n]\setminus [\lceil s\rceil-1]$.
 Clearly,  $\omega$ is a fractional vertex cover of $H_k(n,s)$.   So
 $\nu'(H_k(n,s))=\mu(H_k(n,s))\leq \lceil s\rceil-1$, and the
 assertion of the lemma holds. \qed

\medskip

We also need  two results concerning a famous conjecture of Erd\H{o}s \cite{Er65} on the
matching number of a $k$-graph; both have a requirement on the
number of vertices.
 The first result is due to Frankl (Theorem 1.1 in 
\cite{Fr13}).

\begin{lemma}[Frankl]\label{Fran13}
Let  $k,s$ be integers with $k\ge 2$ and $s\ge 1$. Then, for
any integer $n$ with
 $n\geq (2k-1)s+k$,
$m_0^s(k,n)= {n\choose k}-{n-s+1\choose k}+1$.

\end{lemma}

The second result is due to Frankl and Kupavskii (Theorem 1 in \cite{FK18}).

\begin{lemma}[Frankl and Kupavskii]\label{FrKu18}
Let  $k$ be an integer with $k\ge 2$.  There exists an
absolute constant $s_0\ge 1$ such that, for any integer $s\ge s_0$ and any
integer $n\geq (5k/3-2/3)s$,  $m_0^s(k,n)=  {n\choose k}-{n-s+1\choose k}+1$.
\end{lemma}

We now state and prove the main result of this section, which
essentially says that $f_d^s(k,n)\le f_0^s(k-d, n-d)$, following the method used by Alon et al. in \cite{AFHRS12}.
Recall that for a hypergraph $H$ and $S\subseteq V(H)$,
$N_H(S)=\{T\subseteq V(H)\setminus S: S\cup T\in E(H)\}$. We also view
$N_H(S)$ as a hypergraph with vertex set $V(H)\setminus S$ and edge-set $N_H(S)$.

\begin{lemma}\label{Erdos-frac}
Let $k,d$ be integers with $k\geq 2$ and $1\leq d\leq k-1$, and let $n$ be
a positive integer and  $s$ be a rational constant with $0<s\leq n/k$. Let $H$ be a $k$-graph on $n$ vertices such that, for every $d$-set $S$,
the $(k-d)$-graph $N_{H}(S)$ has a
fractional matching of size at least $s$. Then $H$ has a fractional
matching of size at least $s$.

\end{lemma}

\pf Let $\omega$ be a fractional vertex cover of $H$ with size
$\mu(H)$, and write $V(H)=\{v_1,\ldots,v_n\}$ such that
\begin{itemize}
\item [(1)] $\omega(v_1)\geq \omega(v_2)\geq \cdots\geq \omega(v_n).$
\end{itemize}
Let $H_{\omega}$ be the $k$-graph with vertex set $V(H)$ and edge set
$$E(H_{\omega})=\left\{e\ :\ e\in {V(H)\choose k}\ \mbox{and }\sum_{v\in e}\omega(v)\geq
1\right\}.$$

Then $\omega$ is also a fractional vertex cover of $H_{\omega}$; so $\mu(H_{\omega})\leq
\mu(H)$. Since  every edge of $H$ is also an edge of
$H_\omega$, we have  $\nu'(H_{\omega})\ge \nu'(H)$. Hence,  $\nu'(H_{\omega})=\mu(H_{\omega})\leq \mu(H)=\nu'(H)\le \nu'(H_{\omega})$. Thus, we have

\begin{itemize}
\item [(2)] $ \nu'(H)=\nu'(H_{\omega}).$
\end{itemize}

Let $S=\{v_{n-d+1},\ldots,v_n\}$.
Then, $|S|=d$. Let $w_0:=\frac{1}{d}\sum_{v\in
  S}\omega(v)$, and define $\omega':V(H_{\omega})\rightarrow [0,1]$ such that
\begin{equation*}
 \omega'(v) = \left\{
  \begin{array}{ll}
   \omega(v) , & \hbox{if $v\in V(H_\omega)\setminus S$;} \\
    w_0, & \hbox{if $v\in S$.}
  \end{array}
\right.
\end{equation*}

 We may assume that $w_0<1/k$. For, otherwise, $\nu'(H)=\mu (H)=\sum_{v\in
   V(H)}\omega(v)\ge n\omega_0\ge n/k\ge s$; so the assertion of the
 lemma holds.

 Let $\omega'':V(H)\rightarrow \mathbb{R}^+\cup \{0\}$ be a function such that
\[
\omega''(v)=\frac{\omega'(v)-w_0}{1-kw_0} \quad \mbox{for all $v\in V(H)$. }
\]
Then $\omega''(v)=0$ for $v\in S$.   Note that $N_{H_{\omega}}(S)$ is a
$(k-d)$-graph with vertex set $V(H_{\omega})\setminus S$ (which has $n-d$ vertices).
For any edge $e\in N_{H_{\omega}}(S)$,  since $\omega$ is also a vertex cover of $H_\omega$ and $e\cup S\in E(H_\omega)$, we have $\sum\limits _{v\in e\cup S}\omega(v)\geq 1$. 
Recall that $\omega(v)=\omega'(v)$ for any $v\in V(H)-S$ and $\omega'(x)=0$ for any $x\in S$. So we have
\begin{align*}
\sum_{v\in e}\omega''(v)&=\sum_{v\in e}\frac{\omega'(v)-w_0}{1-kw_0}\\
&=\frac{\sum_{v\in e}\omega'(v)-kw_0}{1-kw_0}\\
&= \frac{\left(\sum\limits _{v\in e\cup S}\omega(v)\right)-kw_0}{1-kw_0}\geq 1.
\end{align*}
Thus, the  function $\omega''$ restricted to $V(H_{\omega})\setminus S$ is  a fractional vertex cover of $N_{H_\omega}(S)$.
Then by  hypothesis and Strong Duality Theorem,  we have
\[
\sum_{v\in V(H_{\omega})\setminus S}\omega''(v)\geq \mu(N_{H_\omega}(S))= \nu'(N_{H_\omega}(S))\geq s.
\]
Recall that $\omega$ is a minimum vertex cover of $H_\omega$. Note that $\nu'(H_{\omega})\leq n/k$; so $k\omega_0\sum_{v\in V(H_{\omega})}\omega(v) \le k\omega_0 (n/k)=n\omega_0$. Hence, we have
\begin{align*}
s\leq \sum_{v\in V(H_\omega)\setminus S}\omega''(v)= \sum_{v\in V(H_\omega)}\omega''(v)&=\frac{\sum_{v\in V(H_\omega)}\omega(v)-nw_0}{1-kw_0}
\leq \sum_{v\in V(H_\omega)}\omega(v)=\nu'(H_{\omega}).
\end{align*}
Thus by (2), $H$ has a fractional matching of size at least $s$.  \qed

\medskip

\section{Proof of Theorem \ref{main-frac}.}

%\textit{Proof of Theorem \ref{main-frac}.}
First, we give a proof of Theorem~\ref{main-frac}.
Let $k,d$ be integers with $k\ge 4$ and $2k/5<d\leq k-1$, and let $n$
be an integer with $n\ge
2k^2$. If $d\ge k/2$ let $s_0=0$, and if $2k/5<d\leq k/2$ let $s_0\ge 1$ be
given as in Lemma~\ref{FrKu18}.  Note that $k-d\geq 2$.  By Lemma~\ref{Erdos-frac},
$f_d^s(k,n)\leq f_0^s(k-d,n-d).$

Since $d<k$,  $\lceil s\rceil \le
(n+k-1)/k<(n-d)/(k-d)$; so  $$f_d^s(k,n)\leq f_0^s(k-d,n-d) \leq m_0^{\lceil
  s\rceil}(k-d,n-d).$$
Therefore, in view of Lemma~\ref{low-bound},
 it suffices to show that $m_0^{\lceil s\rceil}(k-d,n-d)\le {n-d\choose k-d}-{(n-d)-\lceil s\rceil+1\choose
  k-d}+1$ for all $s$ with $s_0<s  \le n/k$ (in which case $\lceil s
\rceil\le  (n-d)/(k-d)$).

We apply Lemma~\ref{Fran13} (when $d\ge k/2$) and Lemma~\ref{FrKu18}
(when $2k/5<d<k/2$) on a $(k-d)$-graph of order $n-d$. Thus,
we need to verify that, for every $s$ with $s_0< s \le n/k$,
$f(d):=(n-d)-[(2(k-d)-1)\lceil s \rceil+(k-d)]\ge 0$ when $d\geq k/2$, and
$g(d):=(n-d)-(5(k-d)/3-2/3)\lceil s\rceil \ge 0$ when $2k/5<d<k/2$.

Note that the first derivatives $f'(d)=2\lceil s\rceil>0$ and
$g'(d)=5\lceil s\rceil /3-1> 0$ when $s>0$. Hence, when $d\ge k/2$,  $f(d)\ge
f(k/2)=n-k-(k-1)\lceil s \rceil \ge 0$,  as $ s \le n/k$ and $n\ge 2k^2$.
When $2k/5<d<k/2$, we have $d\ge (2k+1)/5$. So $g(d)\ge
g((2k+1)/5)=n-(2k+1)/5-(k-1)\lceil s\rceil \ge 0$, as $s 
\le n/k$ and $n\ge 2k^2$.  \qed

\section{Concluding remarks}

R\"odl,  Ruci\'nski, and  Szemer\'edi
\cite{RRE06} determined $f_{k-1}^s (k, n)$ for  $0<s\le n/k$.
For  the entire range $1\le d\le k-2$, K\"{u}hn, Osthus, and Townsend
\cite{KOT14}   proved the following asymptotic result.

\begin{theorem}[Kuhn, Osthus, and Townsend]\label{kot}
Let $k,d$ be integers with $k\geq 3$ and $1\leq d \leq k-2$, and let $0 \leq a \leq \min\{1/(2(k-d)), 1/k\}$.
Then, for positive integers $n$,
\[
f^{an}_d (k, n) \sim\left(1- (1-a)^{k-d} \right){n-d\choose k-d}.
\]
\end{theorem}
Thus,  $f_d^s (k, n)$ is  asymptotically determined when  $1\leq d \leq k-2$
and $s\le n/(2(k-d))$, and when $d\ge k/2$ and $s\in (0, n/k]$.
Theorem~\ref{main-frac} determines $f_d^{s}(k,n)$ exactly
when  $d>2k/5$ and $n,s$ large enough.

\medskip

 For matchings, K\"uhn, Osthus, and Townsend \cite{KOT14} proposed the following conjecture.
\begin{conjecture}[K\"uhn, Osthus, and Townsend]\label{KOT-Conj}
For all $\varepsilon>0$ and all integers $n,k,d,s$ with $1\leq d\leq k-1$ and $1\leq s\leq (1-\varepsilon)n/k$,
\[
m_d^s(k,n)\sim\left(1-(1-s/n)^{k-d}\right){n-d\choose k-d}.
\]

\end{conjecture}

K\"uhn, Osthus, and Townsend \cite{KOT14} proved that Conjecture
\ref{KOT-Conj} holds for $k/2\leq d\leq k-1$.
Han \cite{Han16} showed that this conjecture holds for
$0.42k<d<k/2$.   Alon et al. \cite{AFHRS12} showed  for any two constants $\alpha, \alpha'$ with $0<\alpha'^{1/r}\ll\alpha< 1/k$, where $r$ is a sufficiently large integer, there exists $n_0$ such that for all $n\geq n_0$,  $m_d^{(1-\alpha)n/k}(k,n)\leq f_d^{(1/k-\alpha+\alpha')n}(k,n)$. By Lemma \ref{low-bound}, we have ${n-d\choose k-d}-{n-d-( n/k-\alpha n)\choose k-d}\leq m_d^{(1/k-\alpha)n}(k, n)$.
 Recall that Alon et al. \cite{AFHRS12} proved $m_d^{n/k}(k, n)\sim
\max\{c^*, 1/2\}{n-d\choose k-d}$, where $f_d^{n/k}(k, n)\sim c^*{n-d\choose k-d}$.
Note that for $k\geq 3$ and $2k/5\leq d\leq k-1$, $1-(1-1/k)^{k-d}<1/2$.
As a consequence of Theorem \ref{main-frac} and another result in Alon  et al. (see Theorem 1.1 in \cite{AFHRS12}), we can derive the
following result.

\begin{coro} Let $k,d$ be integers such that $k\ge 2$ and $d>2k/5$. For any constant  $\alpha$  with $0<\alpha < 1/k$,
there exists $n_0$ such that for any $n\geq n_0$,
\[
m_d^{(1/k-\alpha)n}(k, n)\sim{n-d\choose k-d}\left(1-(1-1/k+\alpha)^{k-d}\right),
\]
and
\[
m_d^{n/k}(k,n)\sim \frac{1}{2}{n-d\choose k-d}
\]
\end{coro}

\bigskip

%{\sc Acknowledgement}. We are grateful to an anonymous referee for
%their extensive and thoughtful comments/suggestions and suggesting a
%short proof of Lemma 2.4, which significantly improved
%the exposition and the quality of this note.


\begin{thebibliography}{99}
\addtolength{\baselineskip}{-1ex}


\bibitem{AFHRS12} N. Alon, P. Frankl, H. Huang, V. R\"odl,
  A. Ruci\'nski, and B. Sudakov, Large matchings in uniform hypergraphs
and the conjectures of Erd\H{o}s and Samuels, \emph{J. Combin. Theory Ser. A}, \textbf{119} (2012), 1200--1215.


\bibitem{Er65}P. Erd\H{o}s,   A problem on independent $r$-tuples, \emph{Ann. Univ. Sci. Budapest, E\"otv\"os Sect. Math.}, \textbf{8} (1965),
93--95.

\bibitem{Han16} J. Han,
Perfect matchings in hypergraphs and the Erd\H{o}s matching conjecture, \emph{SIAM J. Discrete Math.,} \textbf{30} (2016), 1351--1357.

\bibitem{Fr13} P. Frankl, Improved bounds for Erd\H{o}s matching conjecture, {\it J. Combin. Theory Ser. A}, {\bf 120} (2013), 1068--1072.


\bibitem{FK18} P. Frankl and A. Kupavskii, The Erd\H{o}s matching conjecture and concentration inequalities, arXiv:1806.08855.

\bibitem{KOT14}
D. K\"uhn, D. Osthus, and T. Townsend, Fractional and integer matchings in uniform hypergraphs, \emph{European J. Combin.},
\textbf{38} (2014), 83--96.


\bibitem{RRE06} V. R\"odl, A. Ruci\'nski, and E. Szemer\'edi, Perfect matchings in
uniform hypergraphs with large minimum degree, \emph{European J. Combin.}, \textbf{27} (2006), 1333--1349.


\bibitem{RRS09} V. R\"odl, A. Ruci\'nski, and E.
Szemer\'edi, Perfect matchings in large uniform hypergraphs with
large minimum collective degree, \emph{J. Comb. Theory Ser. A},
\textbf{116} (2009),  613--636.

\bibitem{TZ12} A. Treglown and Y. Zhao, Exact minimum degree thresholds for perfect matchings in uniform hypergraphs I,
\emph{J. Comb. Theory Ser. A}, \textbf{119} (2012),  1500--1522.

\bibitem{TZ13} A. Treglown and Y. Zhao, Exact minimum degree thresholds for perfect matchings in uniform hypergraphs II,
\emph{J. Comb. Theory Ser. A}, \textbf{120} (2013),  1463--1482.


\end{thebibliography}
\end{document}